 \documentclass[12pt]{article}

\usepackage[T1]{fontenc}

\usepackage{amssymb,latexsym,amsmath}
\usepackage{mathrsfs}
\usepackage{dsfont}
\usepackage{enumitem}
\usepackage{caption}
\usepackage{comment}
\usepackage{bbm}
\usepackage{stix}
\usepackage{appendix}
\usepackage{xcolor}

\usepackage{inputenc}
 \usepackage{makeidx}

\allowdisplaybreaks

\newtheorem{lemma}{Lemma}[section]

\newtheorem{theorem}{Theorem}[section]
\newtheorem{proposition}{Proposition}[section]
\newtheorem{assumption}{Assumption}[section]
\newtheorem{example}{Example}[section]
\newtheorem{remark}{Remark}[section]
\newtheorem{definition}{Definition}[section]

\allowdisplaybreaks

\usepackage{epsfig} 
\usepackage{verbatim}

\def\qed{\hfill $\diamond$}

\newcommand{\sy}[1]{{\color{black} #1}}

\usepackage{amsfonts}
\usepackage{float}
\usepackage{multirow}

\allowdisplaybreaks

\begin{document}
\sloppy
\title{Non-Sequential Decentralized Stochastic Control Revisited: Causality and Static Reducibility}
\author{Omar Mrani-Zentar, Ryan Simpson and Serdar Y\"uksel
\thanks{The authors are with the Department of Mathematics and Statistics,
     Queen's University, Kingston, ON, Canada.}
     }
\maketitle

\begin{abstract}
In decentralized stochastic control (or stochastic team theory) and game theory, if there is a pre-defined order in a system in which agents act, the system is called \textit{sequential}, otherwise it is non-sequential. Much of the literature on stochastic control theory, such as studies on the existence analysis, approximation methods, and on dynamic programming or other analytical or learning theoretic methods, have focused on sequential systems. Many complex practical systems, however, are non-sequential where the order of agents acting is random, and dependent on the realization of solution paths and prior actions taken. The study of such systems is particularly challenging as tools applicable for sequential models are not directly applicable. In this paper, we will first revisit the notion of Causality (a definition due to Witsenhausen and which has been refined by Andersland and Tekenetzis), and provide an alternative representation using imaginary agents. We show that Causality is equivalent to Causal Implementability (and Dead-Lock Freeness), thus, generalizing previous results. We show that Causality, under an absolute continuity condition, allows for an equivalent static model whose reduction is policy-independent. Since the static reduction method for sequential control problems (via change of measures or other techniques), has been shown to be very effective in arriving at existence, structural, approximation and learning theoretic results, our analysis facilitates much of the stochastic analysis available for sequential systems to also be applicable for a class of non-sequential systems. 
\end{abstract}


\section{Introduction}

An increasingly important class of optimal stochastic control problems involve setups where a number of decision makers (DMs) / controllers / agents, who have access to different and local information, are present. Such a collection of decision makers who wish to minimize a common cost function and who has an agreement on the system (that is, the probability space on which the system is defined, and the policy and action spaces) is said to be a {\em stochastic team}. Such problems are also called {\it decentralized stochastic control} problems.

For such problems, classical methods involving single-agent stochastic control often are not applicable due to information structure constraints (and notably due to the fact that the information available to agents acting one after another do not necessarily increase, i.e., there is no perfect recall). Such team problems entail a collection of DMs acting together to optimize a common cost function, but not necessarily sharing all the available information. At each time stage, each DM has only partial access to the global information, which is characterized by the \textit{information structure} of the problem \cite{WitsenhausenSIAM71,wit75}. 

If there is a pre-defined order in which the DMs act, then the team is called a \textit{sequential team}. If each DM's information depends only on primitive random variables, the team is \textit{static}. If at least one DM's information is affected by an action of another DM, the team is said to be \textit{dynamic}. Information structures can be further categorized as \textit{classical}, \textit{partially nested} (or quasi-classical), and \textit{nonclassical}. An information structure is \textit{classical} if the information of decision maker $i$ (DM$^i$) includes all of the information available to DM$^k$ for $k < i$. An information structure is \textit{partially nested}, if whenever the action of DM$^k$, for some $k < i$, affects the information of DM$^i$, then the information of DM$^i$ includes the information of DM$^k$. An information structure that is not partially nested is \textit{nonclassical}. 

Despite extensive studies for sequential decentralized stochastic control problems, which also include dynamic programming methods
(see e.g. under various formulations \cite{NayyarMahajanTeneketzis,WitsenStandard,YukselWitsenStandardArXiv}), and while ubiquitous in many complex engineering systems or other application areas such as economics, research on non-sequential stochastic control has been rather scarce. One could notice that since the foundational works by Witsenhausen  \cite{WitsenhausenSIAM71} and later by Andersland and Teneketzis \cite{AnderslandTeneketzisI}, \cite{AnderslandTeneketzisII} and Teneketzis \cite{Teneketzis2}, there has been modest progress in this area with the exception of several important recent studies such as \cite{alos2008trees,heymann2020kuhn,kadnikov2020jeux} primarily in the game theory community which consider models which either possess perfect recall or extensive tree model representations. 

Static teams are relatively simpler to study: For teams with finitely many DMs, Marschak \cite{mar55} has studied static teams and Radner \cite{Radner} has established connections between person-by-person (pbp) optimality (ie. Nash equilibrium), stationarity, and global optimality. Radner's results were generalized in \cite{KraMar82} by relaxing optimality conditions. The essence of these results is that in the context of static team problems, convexity of the cost function, subject to minor regularity conditions, suffices for the global optimality of pbp optimal solutions. In the particular case of LQG (Linear Quadratic Gaussian) static teams, this result leads to optimality of linear policies \cite{Radner}. Optimality of linear policies also holds for dynamic LQG teams with partially nested information structures through a transformation of the dynamic team to a static one \cite{HoChu}: in \cite{HoChu} dynamic LQG teams with partially nested information structures and in \cite{ho1973equivalence} general dynamic teams with partially nested information structures, satisfying an invertibility assumption, have been shown to be reducible to static team problems, where the aforementioned results for static teams can be applied. Such transformation of dynamic teams to static teams is called {\it{static reduction}}. In the static reduction presented in \cite{HoChu, ho1973equivalence}, given the policies of the DMs, there is a bijection between observations as a function of precedent actions of DMs and  the primitive random variables, and observations generated under the transformations, where now they are only functions of primitive random variables. 

In \cite{wit88}, Witsenhausen introduced a static reduction method for dynamic teams, where observation kernels satisfy an absolute continuity condition. In this static reduction, the probabilistic nature of the problem has been transformed to the cost function by changing the measures of the observations to fixed probability measures. Witsenhausen's static reduction procedure is independent of the policies that precedent DMs choose, and hence, we refer to this type of static reduction as a {\it{policy-independent static reduction}} \cite{sanjari2021optimality}. The policy-independent static reduction is a change of measure argument, a prominent version of which is known as Girsanov's transformation \cite{girsanov1960transforming, benevs1971existence}.

 Since Witsenhausen’s paper \cite{wit88}, the static reduction method for sequential control problems has been shown to be very effective in arriving at existence, structural and approximation results. For existence results building on this approach, we refer the reader to \cite{gupta2014existence,SaldiArXiv2017,YukselWitsenStandardArXiv}, for a dynamic programming formulation to \cite{WitsenStandard} for countable spaces and \cite{YukselWitsenStandardArXiv} for general spaces, for rigorous approximations with finite models to \cite{saldiyuksellinder2017finiteTeam}. 
   
In this paper, in the context of non-sequential decentralized stochastic control, we will first revisit the notion of Causality (a definition due to Witsenhausen and which has been refined by Andersland and Tekenetzis), and provide an alternative representation using imaginary agents. Through this representation, we show that Causality is equivalent to Causal Implementability (and Dead-Lock Freeness), thus, generalizing previous results by Witsenhausen  \cite{WitsenhausenSIAM71} and later by Andersland and Teneketzis \cite{AnderslandTeneketzisI}. Moreover, we will show that under our alternative definition one can obtain a concise proof of the fact that causality leads to solvability measurability \cite{WitsenhausenSIAM71} which is a concept that relates to the well-posedness of the problem. We then study static reducibility for non-sequential stochastic control problems. During our derivations and results, we also present new interpretations for some of the classical results put forward by Witsenhausen, Andersland, and Teneketzis- in terms of the imaginary DM model. It is our hope that research in this field may be revived in the community. 

Static reduction also facilitates stochastic optimality analysis via numerical methods and convex relaxations \cite{saldiyukselGeoInfoStructure} and learning theoretic methods in decentralized stochastic control \cite{ArslanYukselTAC16}. Notably, for non-sequential systems it is particularly challenging to apply dynamic programming or reason via backwards induction. One benefit of static reduction is that even if the cost function associated with a static formulation is not written explicitly or known, the fact that there exists one such cost function, and the fact that the information is static allows one to apply learning algorithms to converge to equilibria \cite{ArslanYukselTAC16}; for continuous spaces this also holds \cite[Section 6.1]{altabaa2023decentralized}. 


With the recognition that the information fields generated by local measurements lead to subtle conditions on solvability and causality, an alternative probabilistic model, based on quantum mechanics, for describing such problems has been proposed by Baras in \cite{Baras1} and \cite{Baras2}. Such quantum team decision theory, for sequential systems, has also been studied in \cite[Section 5.2]{saldiyukselGeoInfoStructure} and \cite{anantharam2025quantum}. In this context, our positive results in the paper under Causality also has positive implications for the applicability of classical probability theory for non-sequential models. 

{\bf Contributions.}
(i) We revisit the notion of Causality (a definition due to Witsenhausen \cite{WitsenhausenSIAM71} and which has been refined by Andersland and Teneketzis \cite{AnderslandTeneketzisI}), and provide an alternative representation using imaginary agents. Through this representation, we show that Causality is equivalent to Causal Implementability (and, therefore, Dead-Lock Freeness) in Theorem \ref{ImaginarySequentialImplicationsEquivalence}, complementing the existing implication relations presented in Theorem \ref{ImplicationsThm} via generalizing the equivalence result of \cite[Theorem 4]{AnderslandTeneketzisI} to general standard Borel models; this addresses an open question on the equivalence between Causality and Causal Implementability for general spaces noted in \cite[p. 251]{Teneketzis2}. Moreover, we provide novel proofs of the fact that Causality implies Dead-lock Freeness and Solvability Measurability (Proposition \ref{C implies DF} and \ref{C implies SM}).

(ii) We show that Causality (C), under an absolute continuity condition, allows for an equivalent static model whose reduction is policy-independent (Theorem \ref{mainThmNSSR}). This result facilitates much of the stochastic analysis available for sequential systems to also be applicable for non-sequential systems. In particular, learning algorithms are guaranteed to converge to optimality for finite space models \cite{ArslanYukselTAC16} and near optimality for Borel models \cite[Section 6.1]{altabaa2023decentralized}, even when the exact form of the reduced model is not explicitly written but that such a reduction is available is known.

\section{Properties of Non-Sequential Decentralized Stochastic Control Problems}\label{SecModel}

Witsenhausen \cite{WitsenhausenSIAM71} defines several properties for non-sequential systems: At a high-level, Witsenhausen defines {\it causality} to be the property that for each $k \leq N$, the events that a DM acting at time $k$ can distinguish are events that can be induced by the decisions of the first $k-1$ DMs and $\omega$. Furthermore, after $k-1$ DMs have acted, the selection of the $k$th DM is determined by $\omega$ and the actions up to time $k-1$. Nonetheless, the DM acting at time $k$ does not need to know that he/she is {\bf DM}$k$. Witsenhausen has proved that causality implies {\it solvability} \cite[Theorem 1]{WitsenhausenSIAM71}, which ensures that all control action variables in the system are well-defined random variables (as proper measurable functions), and Andersland and Teneketzis later introduced \textit{causal-implementability} and \textit{deadlock-freeness} \cite{AnderslandTeneketzisI} properties which pertain to the information available to the DMs, and do not impose restrictions on the ordering of DMs. All of these properties are formally defined and studied in this section. 

Before proceeding with the definitions, we introduce a variation of the intrinsic model presented in \cite{wit75} and some notation that will be used in the definitions.

\subsection{An intrinsic model for non-sequential stochastic control} \label{intrinsic}

The general intrinsic model that we will introduce will be a minor modification of Witsenhausen's sequential model \cite{wit75}. The model consists of:

\begin{enumerate}
\item The number of DMs in the system, $N \in \mathbb{N}$.
\item A measurable space $(\Omega, {\cal F})$ characterizing the underlying event space.
 \item $(\mathbb{U}^{k},{\cal U}^{k})$ denotes the standard Borel space from which $u^{k}$, the $k^{th}$ control action is selected. The measurable space containing the collective action ${\bf u}:= (u^{1}, u^{2}, \cdots, u^{N})$ is denoted by $(\prod_{i=1}^{N}\mathbb{U}^{i}, \otimes_{i=1}^{N} {\cal U}^{i})$.
\item $(\mathbb{Y}^{k}, {\cal Y}^{k})$ denotes the standard Borel space from which $y^{k}$, the measurement {\bf DM}$k$ , takes values from
\item Measurable measurement functions, $y^{i} = \eta^i(\omega, {\bf u})$, which output the measurements of each {\bf DM}$i$, for $i = 1, 2, \cdots, N$. Let $\mathscr{J}^{i}$ denote the sigma field induced on $\Omega \times (\prod_{k=1}^N \mathbb{U}^{k})$ by the measurement function $\eta^i$. 

\item A design constraint which restricts the $N$-tuples of control policies to be measurable given the information structure: $\underline{\gamma} = \{\gamma^{1}, \gamma^{2}, \cdots, \gamma^{N}\}$, with
$\gamma^{k}:(\mathbb{Y}^{k}, {\cal Y}^{k}) \rightarrow (\mathbb{U}^{k}, {\cal U}^{k})$, $k = 1, 2, \cdots, N$, so that $u^k = \gamma^k(y^k)$. We let $\Gamma^k$ denote the set of all such policies for DM$k$, and ${\bf \Gamma} := \prod_{k=1}^N \Gamma^k$ denote the collection of such policies.

\item A probability measure $P$ on $(\Omega, {\cal F})$.
\end{enumerate}

We will also find it useful to refine the exogenous randomness (represented canonically with $\omega \in \Omega$) in terms of the collection $\{\Omega_{s_0},\Omega_{0},\cdots,\Omega_{N}\}$ denoting measurable spaces from which the collection of random variables $\omega_{s_{0}}, \omega_{0}, \omega_{1}, \cdots, \omega_{N}$, all measurable on $(\Omega, {\cal F})$, take values from. Here, $\Omega_{0}$-valued $\omega_{0}$ and $\Omega_{s_0}$-valued $\omega_{s_{0}}$ are the cost and order relevant variables, where $\omega_{0}$ is the cost-relevant uncertainty, and $\omega_{s_{0}}$ is the uncertainty related to the possibly random ordering of DMs.
Let $S_{k}$ denote the set of $k$-DM orderings of $k = 1, 2, ..., N$ and let $\psi : \Omega \times (\prod_{i=1}^N \mathbb{U}^i) \rightarrow S_{N}$ denote an ordering function. Then, $\sigma(\omega_{s_0},u^1,\cdots,u^N)$ is the smallest $\sigma$-field over $\Omega \times_{i=1}^N \mathbb{U}^i$ on which $\psi$ is measurable. If $c: (\omega,u^1,\cdots,u^N) \mapsto c(\omega,u^1,\cdots,u^N)$ is the cost function, $\sigma(\omega_0,\omega_{s_0}, u^1,\cdots,u^N)$ is the smallest $\sigma$-field over $\Omega \times_{i=1}^N \mathbb{U}^i$ on which $c$ is measurable. We note that $(\omega_{0}, \omega_{s_{0}})$ may not be independent, and thus only their joint probability distribution is referred to. The variables $\omega_{1}, \omega_{2}, \cdots, \omega_{N}$ represent the noise in the measurements $y^{i}, i = 1, 2, \cdots, N$.

Given a cost function $c: \Omega_0 \times \Omega_{s_0} \times \prod_{k=1}^N \mathbb{U}^i \to \mathbb{R}_+$, a solution to a stochastic team problem is to identify a policy $\underline{\gamma}$ that achieves $\inf_{\underline{\gamma} \in {\bf \Gamma}} J(\underline{\gamma})$ with
\begin{align}\label{costDefinition}
J(\underline{\gamma}):= E\bigg[c\bigg(\omega_{0}, \omega_{s_{0}},\gamma^1(y^1),\cdots,\gamma^N(y^N)\bigg)\bigg] =E^{\underline{\gamma}}[c(\omega_{0}, \omega_{s_{0}}, {\bf u})].
\end{align}

In the above, $E^{\underline{\gamma}}$ is the expectation over the random variables given policies; likewise, we define $P^{\underline{\gamma}}$ as the probability measure on the system variables given the policy $\underline{\gamma}$, provided that this expectation is well-defined\footnote{As we will see, such measurability properties may not be taken for granted for non-sequential systems}.

\subsection{A generalized classification of information structures for non-sequential systems}

Some of the tediousness that comes with studying non-sequential decentralized stochastic control problems is due to the fact that the information structures of decision makers are not as easily characterized as in sequential problems. Since the indexing of DMs does not represent the order in which they act when studying non-sequential problems, it is necessary to adjust the standard classifications from sequential problems, introduced by Witsenhausen, so that they may also apply to non-sequential problems. In the following, and in the rest of the paper, we denote the ordering of DMs for a particular realized solution path as $s := (s_{1}, s_{2}, ..., s_{N}) \in S_N$, where $S_N$ denotes the set of all possible $N$-agent orderings, and $s_{k}  = i$ means that $DM^{i}$ is the $k$-th DM to act. We offer the definitions of such generalized classifications for non-sequential systems below:

\begin{definition} \textbf{Classifications of Information Structures}
\newline
(i) An information structure is classical if the measurement $y^{s_{i}}$ of the $i^{th}$ $DM$ to act contains all of the information available to $DM^{s_{k}}$ for $k < i$, for $i = 1, 2, ..., N$, and for all possible orderings.
\newline
(ii) An information structure is quasi-classical or partially nested, if whenever $u^{s_{k}}$, for some $k < i$ affects $y^{s_{i}}$ through the measurement function $g_{s_{i}}$, $y^{s_{i}}$ contains all of the information available to $DM^{s_{k}}$. That is $\sigma(y^{s_{k}}) \subset \sigma(y^{s_{i}})$, for all $i = 1, 2, ..., N$, and for all possible orderings. Where $\sigma(y^{s_{k}})$ and  $\sigma(y^{s_{i}})$ denote the smallest sigma algebras such that $y^{s_{k}}$ and $y^{s_{i}}$ are measurable, respectively. 
\newline
(iii) An information structure which is not partially nested is non-classical.

\end{definition}

\subsection{Causality, causal implementability and solvability}
We begin by introducing some notation that will be used in the definitions that follow. Define $\Bar{\mathscr{J}}^{k}$ as the smallest sigma algebra on $\Omega \times \prod_{i=1}^{N}\mathbb{U}^{i}$ such that 
\begin{equation} \label{MeasI}
    \Bar{y}^{k}=\sum_{i=1}^{N}\eta^{i}(\omega,\mathbf{u})1_{\{s_{k}=i\}}
\end{equation}
 is measurable, where we recall $s_{k}$ denotes the identity of the $k^{th}$ DM to act. Later on, $\Bar{y}^{k}$ will serve as the measurement of the $k^{th}$ DM in our imaginary sequential model. We note that the notation above means that $\Bar{y}^{k}=\eta^{i}(\omega,\mathbf{u})$ if $s_{k}=i$, for $i=1,\cdots,N$. Moreover, we note that $\Bar{\mathscr{J}}^{k}$ is different from $\mathscr{J}^{k}$ which is the information available to the DM whose identity is given by $k$. Similarly, one can define \begin{equation} \label{ActionIU}
    \Bar{u}^{k}=\sum_{i=1}^{N}u^{i}1_{\{s_{k}=i\}}.
\end{equation} We denote by $\mathscr{F}^{*}(k)$ the smallest sigma algebra on $\Omega \times \prod_{i=1}^{N}\mathbb{U}^{i}$ such that $(\omega,\Bar{u}^{1},\cdots,\Bar{u}^{k})$ is measurable. Let $S_{k}$ denote the set of $k$-DM orderings of $k = 1, 2, ..., N$. For any $s \in S_{k}, k = 1, 2, ..., N,$ with $s=(s_{1},\cdots,s_{k})$, let $T_{j}^{k}:S_{k} \rightarrow S_{j}$ denote a truncation map that returns the ordering of the first $j$ agents of a $k$-agent ordering, for $j \leq k$.

\subsubsection{Causality (C)}
We begin with causality, which will be of particular importance. We note that definition below is more general than the one provided in \cite{Teneketzis2,WitsenhausenSIAM71}. In particular for the definition presented in \cite{Teneketzis2,WitsenhausenSIAM71} the sigma algebras generated by the measurement functions of all agents have to be trivial whenever the ordering is deterministic. Definition \ref{Alt def causality} only requires that the sigma algebra generated by the measurement variable of the $k^{th}$ agent to act \ref{MeasI} be mathematically determined by $\omega$ as well as the actions of agents that acted before it. Moreover, definition \ref{Alt def causality} requires that the identity of the $k^{th}$ agent to act be mathematically determined by $\omega$ as well as the actions of agents that acted before it. We note, though, that the $k^{th}$ agent to act need not now the order in which they act.

\begin{definition} \label{Alt def causality}
    An information structure possesses Property C if there exists at least one (ordering) map $\psi: \Omega \times \prod_{i=1}^{N} \mathbb{U}^i \rightarrow S_{N}$ such that for any $s = (s_{1}, s_{2}, ..., s_{k}) \in S_{k}, k = 1, 2, ..., N$,

\begin{equation}\label{orderingEqnMap}
    \Bar{\mathscr{J}}^{k} \cap [T_{k}^{N} \circ \psi]^{-1}(s):=\{E\cap [T_{k}^{N} \circ \psi]^{-1}(s)| E\in \Bar{\mathscr{J}}^{k}\} \subseteq \mathscr{F}^{*}(k-1)
\end{equation}

\end{definition}

The term $[T_{k}^{N} \circ \psi]^{-1}(s)$ denotes the preimage of the mapping $T_{k}^{N} \circ \psi : \Omega \times \prod_{i=1}^{N} \mathbb{U}^i  \rightarrow S_{k}$. Therefore, given the ordering $s \in S_{k}$ of the first $k$ DMs, $[T_{k}^{N} \circ \psi]^{-1}(s)$ will be the set of intrinsic outcomes $(\omega, {\bf u}) \in \Omega \times \prod_{i=1}^{N} \mathbb{U}^i $ that maps, through $\psi$, to an $N$-DM ordering in which the ordering of the first $k$ agents is $s$. In particular, $\big\{[T_{k}^{N} \circ \psi]^{-1}(s)|s\in S_{k}\big\}$ is a partition of $\Omega \times \prod_{i=1}^{N} \mathbb{U}^i$ into sets of intrinsic outcomes that lead to the same ordering of the first $k$ agents. 

\begin{definition}
    Consider two sigma algebras $\mathcal{F}_{1}$, $\mathcal{F}_{2}$ which are sets of subsets of some state space $\Omega$. Then, the join of $\mathcal{F}_{1}$ and $\mathcal{F}_{2}$ denoted $\mathcal{F}_{1}\lor \mathcal{F}_{2}$ is the coarsest sigma algebra which contains both $\mathcal{F}_{1}$ and $\mathcal{F}_{2}$ and is given by \[\mathcal{F}_{1} \lor \mathcal{F}_{2}:=\sigma(\{A\cap B\mid A\in\mathcal{F}_{1}\text{, }B\in\mathcal{F}_{2}\})\] 
\end{definition}

\begin{proposition} \label{prop C equivalence}
    Condition (\ref{orderingEqnMap}) is equivalent to 
    \[\mathscr{\Bar{J}}^{k}\lor \sigma(\{[T_{k}^{N} \circ \psi]^{-1}(s)|s\in S_{k}\})\subseteq \mathscr{F}^{*}(k-1),\]
where we use the notation $\lor$ to denote the meet of two sigma algebras, i.e., the coarsest (smallest) sigma algebra containing the two sigma algebras. 
\end{proposition}

\textbf{Proof of Proposition \ref{prop C equivalence}.}  The following lemmas (whose proofs are routine) will be useful.
\begin{lemma} \label{Lemma I}
    Let $(\Omega, \mathcal{F})$ be a measurable space and let $\mathcal{C}$ be a collection of subsets of $\Omega$. Then, $\mathcal{C}\subseteq \mathcal{F}$ if and only if $\sigma(\mathcal{C})\subseteq \mathcal{F} $.
\end{lemma}
\begin{lemma} \label{Lemma II}
     Let $(\Omega, \mathcal{F})$ be a measurable space and let $\mathcal{C}$ be a collection of subsets of $\Omega$ such that $\Omega \in \mathcal{C}$. Then, $\sigma(\{A\cap B\mid A\in\mathcal{F}\text{, }B\in\mathcal{C}\})=\mathcal{F} \lor\sigma(\mathcal{C})$.
\end{lemma}

Now, we have that condition \ref{orderingEqnMap} is equivalent to 
\begin{equation*}
    \{E\cap [T_{k}^{N} \circ \psi]^{-1}(s)| E\in \mathscr{\Bar{J}}^{k}, \text{ } s\in S_{k}\} \subseteq \mathscr{F}^{*}(k-1)
\end{equation*}
By Lemma \ref{Lemma I}, this is equivalent to the condition that 
\begin{equation*}
   \sigma(\{E\cap [T_{k}^{N} \circ \psi]^{-1}(s)| E\in \mathscr{\Bar{J}}^{k}, \text{ } s\in S_{k}\}) \subseteq \mathscr{F}^{*}(k-1)
\end{equation*}
which in turn by Lemma \ref{Lemma II} is equivalent to \[\mathscr{\Bar{J}}^{k}\lor \sigma(\{[T_{k}^{N} \circ \psi]^{-1}(s)|s\in S_{k}\})\subseteq \mathscr{F}^{*}(k-1)\]
\qed

\subsubsection{Equivalence to a sequential model with imaginary decision makers under property C} \label{imaginary model}

In what follows, we view the information structure in terms of {\it imaginary} DMs which are indexed by their time of action rather than their identity. In the following we denote by $\Bar{DM}^{k}$ as the $k^{th}$ DM to act, regardless of the realized solution path. Alternatively, we can say that for any realized solution path with ordering $s = (s_{1}, \cdots, s_{N})$, we have that $DM^{i} = \Bar{DM}^{k}$, when $s_{k}=i$. In this sense, we also change the actions and measurements of each DM to be represented in terms of the order of the actions, rather than the identity of the DM. 

Consider a non-sequential decentralized stochastic control problem described by the intrinsic model in Section \ref{intrinsic}, and suppose that there exists an ordering function $\psi$ such that the information structure possesses property C. Then there exists an equivalent sequential model described by the following:

 \begin{enumerate}
     \item The number of DMs in the system, $N \in \mathbb{N}$.
     \item A measurable space $(\Omega, {\cal F})$ characterizing the underlying event space.
    \item The measurements $\Bar{y}^{k}$ as defined by equation \ref{MeasI}. Let  $(\Bar{\mathbb{Y}}^{k}, \Bar{{\cal Y}}^{k})$ denote the standard Borel space from which $\Bar{y}^{k}$, the measurement of $\Bar{DM}k$ , takes values from. The measurement sets $\Bar{\mathbb{Y}}^{k}$ are related to the original measurement sets $\mathbb{Y}^{k}$ by:
\begin{align}
\label{measSpace}\Bar{\mathbb{Y}}^{k} = \cup_{i=1}^N \mathbb{Y}^{i}. 
\end{align}

\item For $k = 1, \cdots, N$, the actions $\Bar{u}^{k}$ are defined by equation \ref{ActionIU} and can also be written as
  \begin{align}   
      \Bar{u}^{k} &=\gamma^{i}(y^{i}): \text{ if }\psi_{k}(\omega_{s_{0}}, \Bar{u}^{1}, \cdots, \Bar{u}^{k-1}) = s_k = i \label{actionIndexMeasurable} \\
       &=: \Bar{\gamma}^{i}(\Bar{y}^{i},s_i) \label{actionIndex}
      \end{align}
      $(\Bar{\mathbb{U}}^{k},\Bar{{\cal U}}^{k})$ denotes the standard Borel space from which $\Bar{u}^{k}$, the $k^{th}$ control action is selected. The measurable space containing the collective action $\Bar{{\bf u}}:= (\Bar{u}^{1}, \Bar{u}^{2}, \cdots, \Bar{u}^{N})$ is denoted by $(\prod_{i=1}^{N}\Bar{\mathbb{U}}^{i}, \otimes_{i=1}^{N} \Bar{{\cal U}}^{i})$. The action sets $\Bar{\mathbb{U}}^{k}$ are related to the original action sets $\Bar{U}^{k}$ by:
      \begin{align}
\label{measSeqImSpaceActSpace}\Bar{\mathbb{U}}^{k} = \cup_{i=1}^N \mathbb{U}^{i},
\end{align}
    
    
    \item A probability measure $P$ on $(\Omega, {\cal F})$.
\end{enumerate}

\begin{remark}
    Let $\psi_{j}:\Omega\times \prod_{i=1}^{N}(\mathbb{U}^{i})\rightarrow \{1,\cdots,N\}$ be a function which maps $(\omega,u^{1},\cdots,u^{N})$ to the identity of the $j^{\text{th}}$ player, i.e.,  $\psi_{j}(\omega,u^{1},\cdots,u^{N}):=s_{j}(\omega,u^{1},\cdots,u^{N})$. Note that when Definition \ref{Alt def causality} holds for some $\psi$ it implies that $\psi_{1}(\omega,u^{1},\cdots,u^{N})=\psi_{1}(\omega)$ and for $j\in \{2,\cdots,N\}$ $\psi_{j}(\omega,u^{1},\cdots,u^{N})=\psi_{j}(\omega,\Bar{u}^{1},\cdots,\Bar{u}^{j-1})$ where $\Bar{u}^{k}$ denotes the action of the $k^{\text{th}}$ DM.
\end{remark}
\begin{remark}
    One can express the original cost function as $c(\omega,\omega_{s_{0}},u^{1},\cdots,u^{N})=\Bar{c}(\omega,\omega_{s_{0}},\Bar{u}^{1},\cdots,\Bar{u}^{N})$ since one can write for all $i=\{1,\cdots,N\}$, $u^{i}=\sum_{k=1}^{N}\Bar{u}^{k}1_{\{s_{k}=i\}}$. At times, it will also be useful to express the cost as function of the ordering and measurements $\Bar{c}(\omega,\omega_{s_{0}},s, \Bar{u}^{1},\cdots,\Bar{u}^{N},\Bar{y}^{1},\cdots,\Bar{y}^{N})$
\end{remark}

We provide an example (which satisfies Property C) in the following.

\begin{example}\label{Ornek1}

We have $N=3$; $\mathbb{U}^{1} = \mathbb{U}^{2} = \mathbb{U}^{3} = \{0, 1\}$; $\mathbb{Y}^{1} = \mathbb{Y}^{2} = \mathbb{Y}^{3} = \{0, \frac{1}{2}, 1\}$, and $\Omega_0 = \Omega_{s_0} = \{0, 1\}$. Additionally,
\begin{enumerate}

 \item[(a)]
 \begin{align}
& \mathscr{J}^{1} = \sigma\bigg\{ \{(\omega, {\bf u}): \omega_{s_{0}} + \omega_{0} = 0\}, \nonumber \{(\omega, {\bf u}): \omega_{s_{0}}\omega_{1}u^{2} = 1\}\bigg\} \nonumber \\
& \mathscr{J}^{2} =\sigma\bigg\{\{(\omega, {\bf u}): \omega_{s_{0}}\omega_{0} = 1\}, \nonumber \{(\omega, {\bf u}): \omega_{s_{0}} + u^{1}u^{3} = 0\} \bigg\}\nonumber \\
& \mathscr{J}^{3} =\sigma\bigg\{ \{(\omega, {\bf u}): \omega_{3}u^{1} = 1\}, \nonumber \{(\omega, {\bf u}): u^{1} + u^{2} = 0\}\bigg\},
\end{align}
where for a collection of sets ${\cal A}$, $\sigma({\cal A})$ denotes the $\sigma$-field generated by the collection over $\Omega \times \mathbb{U}^1 \times \mathbb{U}^2 \times \mathbb{U}^2$. The above are induced by the following measurement functions:

$$y^{1} = \eta^{1}(\omega_{0}, \omega_{s_{0}}, \omega_{1}, u^{2}, u^{3}) = \begin{cases}
1, & \omega_{s_{0}} + \omega_{0} = 0 \\
0, & \omega_{s_{0}}\omega_{1}u^{2} = 1 \\
1/2, & otherwise
\end{cases},$$ 
$$y^{2} = \eta^{2}(\omega_{0}, \omega_{s_{0}}, \omega_{2}, u^{1}, u^{3}) = \begin{cases}
1, & \omega_{s_{0}}\omega_{0} = 1 \\
0, & \omega_{s_{0}} + u^{1}u^{3} = 0  \\
1/2, & otherwise
\end{cases},$$ 
$$y^{3} = \eta^{3}(\omega_{0}, \omega_{s_{0}}, \omega_{3}, u^{1}, u^{2}) = \begin{cases}
1, & \omega_{3}u^{1} = 1 \\
0, & u^{1} + u^{2} = 0 \\
1/2, & otherwise
\end{cases}$$

 \item[(b)] The ordering function
 $\psi := (\psi_{1}, \psi_{2}, \psi_{3}),$ is such that: \\
 $s_1=\psi_{1}(\omega_{s_{0}}) = \begin{cases}
    1, & \omega_{s_{0}} = 0 \\
    2, & \omega_{s_{0}} = 1\\
    \end{cases}$,   $s_2=\psi_{2}(\omega_{s_{0}}, u^{s_{1}}) = 
    \begin{cases}
    1, & \omega_{s_{0}} = 1 \\
    3, & \omega_{s_{0}} = 0, u^{s_{1}} = 1 \\
    2, & otherwise
    \end{cases},$
    $$s_3=\psi_{3}(\omega_{s_{0}}, u^{s_{1}}, u^{s_{2}}) = \begin{cases}
    2, & \omega_{s_{0}} = 0, u^{s_{1}} = 1 \\
    3, & otherwise
    \end{cases}$$
\end{enumerate}

 \item[(c)]
We have, with the $\eta^1, \eta^2, \eta^3$ defined above\\
$\Bar{y}^{1} = y^{s_1} = \eta^{1}(\omega_{0}, \omega_{s_{0}}, \omega_{1}, u^{2}, u^{3}) 1_{\{\omega_{s_{0}} = 0\}} + \eta^{2}(\omega_{0}, \omega_{s_{0}}, \omega_{2}, u^{1}, u^{3})1_{\{\omega_{s_{0}} = 1\}}$,   

$\Bar{y}^{2} = y^{s_2} = \eta^{1}(\omega_{0}, \omega_{s_{0}}, \omega_{1}, u^{2}, u^{3})1_{\{\omega_{s_{0}} = 1\}}  +\eta^{2}(\omega_{0}, \omega_{s_{0}}, \omega_{2}, u^{1}, u^{3}) (1 - 1_{\{\omega_{s_{0}} = 0, u^{s_{1}} = 1\}} - 1_{\{\omega_{s_{0}} = 1\}}) + \eta^{3}(\omega_{0}, \omega_{s_{0}}, \omega_{3}, u^{1}, u^{2})1_{\{\omega_{s_{0}} = 0, u^{s_{1}} = 1\}}$,

$\Bar{y}^{3} = y^{s_3} = \eta^{2}(\omega_{0}, \omega_{s_{0}}, \omega_{2}, u^{1}, u^{3}) 1_{\{\omega_{s_{0}} = 0, u^{s_{1}} = 1\}} + \eta^{3}(\omega_{0}, \omega_{s_{0}}, \omega_{3}, u^{1}, u^{2}) (1-1_{\{\omega_{s_{0}} = 0, u^{s_{1}} = 1\}})$


Finally, the actions are generated in accordance with (\ref{ActionIU}) so that:
  \begin{align}\label{ActionIUOrnek}
      \Bar{u}^{1} &=  \gamma^1(y^1) 1_{\{\omega_{s_{0}} = 0\}}  + \gamma^2(y^2)1_{\{\omega_{s_{0}} = 1\}} \nonumber \\
      \Bar{u}^{2} &=  \gamma^1(y^1) 1_{\{\omega_{s_{0}} = 1\}}  + \gamma^2(y^2)(1 - 1_{\{\omega_{s_{0}} = 0, \Bar{u}^{1}= 1\}}-1_{\{\omega_{s_{0}}=1\}}) + \gamma^3(y^3)1_{\{\omega_{s_{0}} = 0, \Bar{u}^{1} = 1\}} \nonumber \\
      \Bar{u}^{3} &=  \gamma^2(y^2) 1_{\{\omega_{s_{0}} = 0, \Bar{u}^{1} = 1\}}  + \gamma^3(y^3)(1-1_{\{\omega_{s_{0}} = 0, \Bar{u}^{1} = 1\}}) \nonumber 
  \end{align}
\qed

    \end{example}

\begin{remark}
It is important to note that under the model presented above it is not necessary for the $k^{th}$ DM to act to know which DMs have acted previously or that they are the $k^{th}$ DM. However, the universe or \textit{umpire} (with sigma-field $\mathscr{F}^{*}(k-1)$) knows the identity of the $k^{th}$ DM. 
\end{remark}
Our interpretation in terms of 'imaginary' agents leads to a sequential model as the ordering of the imaginary agents is deterministic. Therefore, we have the following:
\begin{proposition}
For a non-sequential system with Property C and ordering function $\psi$ (\ref{orderingEqnMap}), for any admissible team policy, there exists a sequential system with measurement and action spaces as in (\ref{measSpace}, \ref{measSeqImSpaceActSpace}), measurement realization given as in (\ref{MeasI}), policies given with (\ref{actionIndexMeasurable}) and (\ref{actionIndex}), and actions in (\ref{ActionIU}).
\end{proposition}

An important feature of Property C is that it implies that after the first $k-1$ agents have taken their actions, the $k^{th}$ agent to act is thereby determined mathematically. The latter is referred to as random causal sequentiality:

 \begin{definition} [Random Causal Sequentiality (RCS)]:
 A non-sequential system with ordering function $\psi$ possesses property RCS if for any $s = (s_{1}, s_{2}, ..., s_{k}) \in S_{k}, k = 1, 2, ..., N$,

\begin{equation}
     [T_{k}^{N} \circ \psi]^{-1}(s) \subseteq \mathscr{F}^{*}(k-1)
\end{equation}

 \end{definition}

\subsubsection{Solvability (SM) and Causal-Implementability (CI)}

In order to ensure that the control actions are well-defined random variables, and thus to define the expected cost function under a given policy in the non-sequential setup, the information structure should possess a property called \textit{solvability-measurability} (SM) \cite{WitsenhausenSIAM71}.
\begin{definition} [Solvability-Measurability (SM)] \cite{WitsenhausenSIAM71}
Property SM holds when, for each $\gamma \in \Gamma$ and $\omega \in \Omega$, there exists one and only one ${\bf u} \in \prod_{i=1}^N \mathbb{U}^i$ satisfying the closed loop equations
$$P_{i}({\bf u}) = \gamma^{i}(y^{i}), i = 1, 2, ..., N \quad \quad \Rightarrow \quad \quad u^{i}= \gamma^{i}(\eta^{i}(\omega_{0}, \omega_{i}, u^{1}, u^{2}, ..., u^{N}))$$
Here $P_{i}$ is the projection (or restriction) of the product $\prod_{k=1}^{N}\mathbb{U}^{k}$ to $\mathbb{U}^{i}$. Furthermore, the Property SM requires that the solution map denoted by $M^{\gamma}: (\Omega, {\cal F}) \rightarrow (\prod_{k=1}^{N}\mathbb{U}^{k},{\cal B}(\prod_{k=1}^{N}\mathbb{U}^{k}))$ be measurable\footnote{We note that this definition of SM is presented differently than in \cite{WitsenhausenSIAM71}, in order to accommodate use of the measurements $y^{i}, i = 1, 2, ..., N$. The closed loop equations in \cite{WitsenhausenSIAM71} are given by: $P_{i}({\bf u}) = \gamma^{i}(\omega, {\bf u}), i = 1, 2, ..., N$.}. 
\end{definition}

Property SM guarantees that the control actions $u^{i}$, $i = 1, 2, ..., N$, are well-defined random variables, given an admissible control policy $\gamma \equiv (\gamma^{1}, \gamma^{2}, ..., \gamma^{N}).$ Property SM also ensures that the expectation of the cost function is well-defined under a given policy.

A further property in non-sequential stochastic control, called \textit{causal-implementability} (CI) defined by Andersland and Teneketzis in \cite[Definition 2]{AnderslandTeneketzisI}, is necessary and sufficient to ensure that the information available to any DM does not depend on the actions of DMs who act after it \cite[Theorem 1]{AnderslandTeneketzisI}:
\begin{definition}

An information structure possesses Property CI when there exists at least one map $\psi: \Omega \times \prod_{i=1}^N \mathbb{U}^i \rightarrow S_{N}$ such that for any $ k \in \{1, 2, ..., N\}$,

\begin{equation} \label{sigmafieldCI}
     \Bar{\mathscr{J}}^{k} \subseteq \mathscr{F}^{*}(k-1)
\end{equation}
\end{definition}

In \cite[Definition 2]{AnderslandTeneketzisI} property CI is defined in terms of the existence of an ordering map $\psi$ such that for any $(\omega,\mathbf{u})$ the following condition holds
\begin{equation} \label{CI-Teneketzis}
    \mathscr{J}^{s_{k}}\bigcap [\mathcal{P}_{T^{N}_{k-1}(s)}]^{-1}\big(\mathcal{P}_{T^{N}_{k-1}(s)}(\omega,\mathbf{u})\big)\subseteq \bigg\{ \varnothing, [\mathcal{P}_{T^{N}_{k-1}(s)}]^{-1}\big(\mathcal{P}_{T^{N}_{k-1}(s)}(\omega,\mathbf{u})\big)\bigg\}
\end{equation}
where $s=\psi(\omega,\mathbf{u})$ and $\mathcal{P}_{s}(\omega,\mathbf{u})=(\omega,u^{s_{1}},\cdots,u^{s_{k}})$, $\mathcal{P}_{\emptyset}(\omega,\mathbf{u})=\omega$. We note that $\mathscr{J}^{s_{k}}$ is not formally introduced in \cite{AnderslandTeneketzisI} which causes some notational ambiguity. We also note that because of the use of singletons in (\ref{CI-Teneketzis}), it is more suitable for countable spaces whereas (\ref{sigmafieldCI}) is more suitable for standard Borel spaces. If the sigma field $\mathscr{J}^{s_{k}} $ is taken as referring to $\Bar{\mathscr{J}}^{k} $, then (\ref{CI-Teneketzis}) implies (\ref{sigmafieldCI}) when all spaces are finite or countable as shown in the following proposition. 
\begin{proposition}
Suppose $\Omega$ and all action space $\mathbb{U}^{1},\cdots,\mathbb{U}^{N}$ are countable. Then, (\ref{CI-Teneketzis}) implies (\ref{sigmafieldCI}). 
\end{proposition}
\textbf{Proof.} Suppose (\ref{CI-Teneketzis}) holds. Then, we have that 
\begin{eqnarray*}
\sigma\bigg(\bar{\mathscr{J}}^{k}\bigcap \big\{ [\mathcal{P}_{T^{N}_{k-1}(s)}]^{-1}\big(\mathcal{P}_{T^{N}_{k-1}(s)}(\omega,\mathbf{u})\big)\mid \omega,\mathbf{u},s\big\}\bigg)\subseteq\sigma \bigg\{ \varnothing, \big\{[\mathcal{P}_{T^{N}_{k-1}(s)}]^{-1}\big(\mathcal{P}_{T^{N}_{k-1}(s)}(\omega,\mathbf{u})\big)\mid \omega, \mathbf{u},s \big\}\bigg\}
\end{eqnarray*}
By Lemma \ref{Lemma II} this implies that 

\begin{eqnarray*}
\bar{\mathscr{J}}^{k}\lor \sigma\big\{ [\mathcal{P}_{T^{N}_{k-1}(s)}]^{-1}\big(\mathcal{P}_{T^{N}_{k-1}(s)}(\omega,\mathbf{u})\big)\mid \omega,\mathbf{u},s \big\}\subseteq \sigma\big\{[\mathcal{P}_{T^{N}_{k-1}(s)}]^{-1}\big(\mathcal{P}_{T^{N}_{k-1}(s)}(\omega,\mathbf{u})\big)\mid \omega, \mathbf{u},s \big\}
\end{eqnarray*}
Thus
\begin{eqnarray*}
    \bar{\mathscr{J}}^{k}\subseteq \sigma\big\{[\mathcal{P}_{T^{N}_{k-1}(s)}]^{-1}\big(\mathcal{P}_{T^{N}_{k-1}(s)}(\omega,\mathbf{u})\big)\mid \omega, \mathbf{u},s \big\}=\mathscr{F}^{*}(k-1)
\end{eqnarray*}
\qed

Following the definition of Property CI, the following property, called \textit{deadlock-freeness} (DF) (also introduced by \cite{AnderslandTeneketzisI}), is a more easily interpretable and equivalent (see \cite{AnderslandTeneketzisI}) {\it high-level} description of Property CI:
\begin{definition}\label{DFDefinition} [Deadlock-Freeness (DF)] \cite{AnderslandTeneketzisI}
An information structure possesses Property DF if for each $\gamma \in \Gamma$, and for every $\omega \in \Omega$, there exists an ordering of $\gamma$'s $N$ control laws, say $\gamma^{s_{1}(\omega)}, \gamma^{s_{2}(\omega)}, ..., \gamma^{s_{N}(\omega)}$, such that no control action $u^{s_{i}(\omega)}, i = 1, 2, ..., N$, depends on the control actions that follow.
\end{definition}


While the next result has already been established \cite{AnderslandTeneketzisI}, due to the brevity of the proof under our equivalent formulation, a proof is included below.

 \begin{proposition} \label{C implies DF}
    Property C (Definition \ref{Alt def causality}) implies DF.
\end{proposition}
    
    \textbf{Proof.} Let $\psi$ be such that (\ref{orderingEqnMap}) holds. Let $\omega\in \Omega$ and $\gamma \in \Gamma$. The proof will proceed by induction. For the base case $n=1$, it follows from Proposition \ref{prop C equivalence} that $s_{1}$ is uniquely determined by $\omega$ and so is the information available to the first DM: $\Bar{y}^{1}$. Thus, $\Bar{u}^{1}=u^{s_{1}(w)}=\gamma^{s_{1}(w)}(\Bar{y}^{1})$ is independent from the control actions that follow. Now, for the inductive step, suppose for some $n\in\{1,\cdots, N-1\}$ we have that the controls $\Bar{u}^{1},\cdots,\Bar{u}^{n}$ are independent of the control actions that follow. Moreover, suppose that $s_{1},\cdots,s_{n}$ are all functions of $\omega$. Again, by Proposition \ref{prop C equivalence} we get that \[\mathscr{\Bar{J}}^{n+1}\lor \sigma(\{[T_{n+1}^{N} \circ \psi]^{-1}(s)|s\in S_{k}\})\subseteq \mathscr{F}^{*}(n)\] In particular, $\sigma(\{[T_{n+1}^{N} \circ \psi]^{-1}(s)|s\in S_{k}\})\subseteq \mathscr{F}^{*}(n)$. Thus, $s_{n+1}=s_{n+1}(s_{1},\cdots,s_{n},\Bar{u}^{1},\cdots,\Bar{u}^{n})=s_{n+1}(\omega,\Bar{u}^{1},\cdots,\Bar{u}^{n})=s_{n+1}(\omega,\gamma^{s_{1}(w)}(\Bar{y}^{1}),\cdots,\gamma^{s_{n}(w)}(\Bar{y}^{n}))=s_{n+1}(\omega)$. Similarly, we have that $\mathscr{\Bar{J}}^{n+1} \subseteq \mathscr{F}^{*}(n)$  and thus $\Bar{u}^{n+1}=u^{s_{n+1}(w)}=\gamma^{s_{n+1}(w)}(\Bar{y}^{n+1})$ is independent from the control actions that follow.
    \qed

\begin{proposition} \label{C implies SM}
    Property C (Definition \ref{Alt def causality}) implies SM.
\end{proposition}
\textbf{Proof.} Let $\psi$ be such that \ref{orderingEqnMap} holds. Let $\omega\in \Omega$ and $\gamma \in \Gamma$. It follows from the proof of Proposition \ref{C implies DF} that there exists one and only one ${\bf u} \in \prod_{i=1}^N \mathbb{U}^i$ satisfying the closed loop equations
$$P_{i}({\bf u}) = \gamma^{i}(y^{i}), i = 1, 2, ..., N \quad \quad \Rightarrow \quad \quad u^{i}= \gamma^{i}(\eta^{i}(\omega_{0}, \omega_{i}, u^{1}, u^{2}, ..., u^{N}))$$ Now, we note that the first action to be taken can be written as 

\begin{eqnarray*}
    \Bar{u}^{1}=\sum_{i=1}^{N} 1_{\{s_{1}(\omega)=i\}}\gamma^{i}\bigg(\sum_{i=1}^{N} 1_{\{s_{1}(\omega)=i\}}\eta^{i}(\omega)\bigg)
\end{eqnarray*}
Thus $\Bar{u}^{1}$ is a measurable function of $\omega, \gamma$. Similarly, one gets that the functions $\Bar{u}^{2},\cdots,\Bar{u}^{N}$ and $s_{1},\cdots, s_{N}$ are measurable. Hence, for all $i\in \{1,\cdots, N\}$ $u^{i}=\sum_{k=1}^{N}\Bar{u}^{k}1_{\{s_{k}=i\}}$ is measurable which implies that the solution map $M^{\gamma}$ is measurable.\qed

 To gain further insight into the definitions above, consider the following two examples from \cite{Teneketzis2}\footnote{In addition to the Properties SM, DF, C, and CI, \cite{AnderslandTeneketzisII} defines analogues of these properties which are policy-dependent in the sense that they apply only for a given policy or a subset of policies but not necessarily for all policies.}. 

 \begin{example}
Let $\Omega=\mathbb{U}^1=\mathbb{U}^2=\mathbb{U}^3=\{0,1\}$ and $\mathscr{F}^{1}=\sigma\{ \omega(1-u^2)u^3\}$, $\mathscr{F}^{2}=\sigma\{\omega(1-u^3)u^1\}$, $\mathscr{F}^{3}= \sigma\{\omega(1-u^1)u^2\}$. Consider $\omega=1$. In this case under admissible policies for each $i=1,2,3$, $u^{i}$ depends on $u^{-i}:=\{u^1,u^2,u^3\} \setminus \{u^i\}$. This system has a deadlock, since no DM can act.
\end{example}

 \begin{example}
Let $\Omega=\mathbb{U}^1=\mathbb{U}^2=\{0,1\}$ and $\mathscr{F}^{1}=2^{\Omega \times \mathbb{U}^2}$, $\mathscr{F}^{2} = 2^{\Omega \times \mathbb{U}^1}$ (where the notation $2^\mathbb{U}$ denotes the power set, that is the collection of all subsets of the set $\mathbb{U}$). Consider the following team policy:
\[\gamma^1(\omega,u^2)= 0 \times 1_{\{u^2=0\}} + 1 \times 1_{\{u^2=1\}}, \quad \quad \gamma^2(\omega,u^1)= 0 \times 1_{\{u^1=0\}} + 1 \times 1_{\{u^1=1\}},\]
 For this design, consider the realization $\omega=0$. In this case, $(\omega,u^1,u^2)=(0,0,0)$ as well as $(0,1,1)$ are acceptable realizations given the policy stated above. A similar setting occurs for $\omega=1$, since $(1,0,0)$ and $(1,1,1)$ are acceptable realizations. Hence, for a given cost function $c$, there does not exist, in general, a well-defined (measurable) cost realization variable $c(\omega,u^1,u^2)$ under this policy, and the expectation $E[c(\omega,u^1,u^2)]$ is not well defined given the policy $(\gamma^1,\gamma^2)$. This system is not {\it solvable}.
\end{example}


\subsection{Equivalence of Property C and Property CI (and Property DF)}

We first present the following, due to Witsenhausen, and Andersland and Teneketzis:

\begin{theorem}\label{ImplicationsThm} \cite[Theorem 1]{WitsenhausenSIAM71}\cite[Theorems 1 and 2]{AnderslandTeneketzisI} The following relationships hold:
\[C \Rightarrow CI \iff DF \Rightarrow SM\]
\end{theorem}

We note that both of the relations in $C \Rightarrow CI \iff DF$ follow essentially directly from the (alternative and equivalent) $\sigma$-field characterizations of these definitions. The relations presented above were first established in \cite{AnderslandTeneketzisI}; Witsenhausen proved in \cite[Theorem 1]{WitsenhausenSIAM71} that $C \Rightarrow SM$, which is proved indirectly by the fact that $C \Rightarrow CI \Rightarrow SM$

Property CI, in view of the definitions and discussions above, states that for any possible ordering the information available to each DM at the time they take their action cannot depend on the action of any DM that has not yet acted. Note that this definition does not impose any restrictions on the information available to select the ordering of DMs, unlike in Property C, which requires that the ordering function has nested information fields. The only difference between Property C and Property CI is indeed that under C, $\psi$ is causal whereas this is not the case under CI. We note that Property SM does not imply Property C.

That Property C and CI are near equivalent under countability conditions has been reported in \cite[Theorem 4]{AnderslandTeneketzisI}. We now present a more general result, which shows that Causality and Properties CI and DF are essentially equivalent; this addresses an open question of Teneketzis on the equivalence between Causality and Causal Implementability for general spaces, noted in \cite[p. 251]{Teneketzis2}.

\begin{theorem}\label{ImaginarySequentialImplicationsEquivalence}
For the case with standard Borel measurement and action spaces, Property CI implies Property C, in the sense that for every non-sequential system with Property CI one can construct a non-sequential system with Property C such that every admissible policy in the original problem is admissible for the new problem and achieves the same cost performance. Therefore, (together with Theorem \ref{ImplicationsThm}) we have
\[C \iff CI \iff DF \Rightarrow SM\]
\end{theorem}

\textbf{Proof.} We have that Property C implies Property CI by Theorem \ref{ImplicationsThm}. We prove the other direction. Here we use the imaginary sequential model. By the Borel isomorphism theorem, one can always construct an equivalent problem such that the spaces $\mathbb{Y}^{1},\cdots, \mathbb{Y}^{N}$ are pairwise disjoint. Let $s\in S_{N}$. Thus, under our construction, by observing $\Bar{y}^{k}$, the umpire can identify the identity of the $k^{th}$ DM to act. We, thus, obtain that 
\begin{eqnarray*}
    \Bar{\mathscr{J}}^{k}\cap\psi_{k}^{-1}(s)\subseteq \Bar{\mathscr{J}}^{k} \subseteq \mathscr{J}^{*}(k-1)
\end{eqnarray*}
for all $k\in\{1,,\cdots,N\}$ where $\psi_{k}(\omega,\mathbf{u})=s_{k}$. Here the last inclusion follows from (\ref{sigmafieldCI}). Because $\mathscr{J}^{*}(0)\subseteq\mathscr{J}^{*}(1)\subseteq\cdots \mathscr{J}^{*}(N)$, we get that for all $j\in \{1,\cdots,N\}$ $\psi_{j}^{-1}(s)\subseteq \mathscr{J}^{*}(k)$ for all $k\in \{j-1,\cdots,N\}$. Thus, for any $k\in\{1,\cdots,N\}$ we get that $\displaystyle [T_{k}^{N} \circ \psi]^{-1}(s)=\cap_{j=1}^{k}\psi_{j}^{-1}(s)\subseteq \mathscr{J}^{*}(k-1)$. The latter in combination with property CI entails that (\ref{orderingEqnMap}) holds. \qed

\section{Static Reduction of Non-Sequential Decentralized Stochastic Control}\label{StaticRedSec}
\subsection{Policy-Independent Static Reduction under Causality}
We now study non-sequential stochastic control problems which admit a policy-independent static reduction. Let us consider an $N$-DM non-sequential decentralized stochastic control problem with measurements $(y^{i})_{i=1, \cdots N}$ and actions $(u^{i})_{i=1, \cdots, N}$, and with its' equivalent imaginary sequential system having measurements $(\Bar{y}^{i})_{i=1, \cdots, N}$ and actions $(\Bar{u}^{i})_{i=1, \cdots, N}$. The following assumption will be useful
\begin{assumption} \label{abs continuity assumption} For all $i\in\{1,\cdots,N\}$
   \begin{align} \label{cont}
       & P(\Bar{y}^{i} \in A|\omega_{s_{0}}, \omega_{0}, \Bar{y}^{1}, \cdots, \Bar{y}^{i-1}, \Bar{u}^{1}, \cdots, \Bar{u}^{i-1}) =\int_{A}f_{i}(\Bar{y}^{i}, \omega_{s_{0}}, \omega_{0}, \Bar{y}^{1}, \cdots, \Bar{y}^{i-1}, \Bar{u}^{1}, \cdots, \Bar{u}^{i-1})\Bar{Q}_{i}(d\Bar{y}^{i}),
       \end{align}
 where $f_{i}$ is the {\it Radon-Nikodym derivative}  of $P(\Bar{y}^{i} \in A|\omega_{s_{0}}, \omega_{0}, \Bar{y}^{1}, \cdots, \Bar{y}^{i-1}, \Bar{u}^{1}, \cdots, \Bar{u}^{i-1}) $ with respect to $\Bar{Q}_{i}(d\Bar{y}^{i})$. \footnote{
As Witsenhausen notes while studying \cite[Eqn (4.2)]{wit88}, when the measurement variables take values from countable set, a reference measure can be constructed (e.g., $Q_t(z) = \sum_{i \geq 1} 2^{-i} 1_{\{z = m_i\}}$ where $\mathbb{Y}^t=\{m_i, i \in \mathbb{N}\}$) so that the absolute continuity condition always holds.}
\end{assumption}
In terms of the original problem, the following proposition offers a sufficient condition for Assumption \ref{abs continuity assumption} to hold.

\begin{proposition}\label{nonSeqReducStandardSpace} Suppose that an $N$-DM non-sequential stochastic control problem is such that that for all $i\in \{1,\cdots,N\}$ $P(dy^{i} |\omega_{s_{0}}, \omega_{0}, y^{1}, \cdots, y^{i-1}, u^{1}, \cdots, u^{i-1})\ll \zeta_{i}(dy^{i})$ where $\zeta_{i}\in \mathcal{P}(\mathbb{Y}^{i})$. Then, Assumption \ref{abs continuity assumption} holds.
\end{proposition}

\textbf{Proof.} It is sufficient to note that $P(\Bar{y}^{i} \in A|\omega_{s_{0}}, \omega_{0}, \Bar{y}^{1}, \cdots, \Bar{y}^{i-1}, \Bar{u}^{1}, \cdots, \Bar{u}^{i-1}) \ll \frac{1}{N}\sum_{i=1}^{N}\zeta_{i}$. \qed

Before we proceed to the main theorem of this section we note that under Property C, the measurement kernels in the imaginary sequential model are policy independent. For the sake of clarity, we state this mathematically in the following Lemma. 

\begin{lemma}\label{ImaginarySequentialModel}
Under any team policy $\gamma$, under Property C
\begin{align}
& P^{\gamma}(d\Bar{y}^{k}|\omega_{0}, \omega_{s_{0}}, \Bar{y}^{1}, \Bar{u}^{1}, \cdots, \Bar{y}^{k-1}, \Bar{u}^{k-1}) = P(d\Bar{y}^{k}|\omega_{0}, \omega_{s_{0}}, \bar{y}^{1}, \Bar{u}^{1}, \cdots, \Bar{y}^{k-1}, \Bar{u}^{k-1})
\end{align}
That is, the conditional probabilities defining measurement kernels are policy independent.
\end{lemma}
\textbf{Proof.}
We have that, for all $1 \leq k \leq N$
\begin{align}
& P^{\gamma}(d\Bar{y}^{k}|\omega_{0}, \omega_{s_{0}}, \Bar{y}^{1}, \Bar{u}^{1}, \cdots, \Bar{y}^{k-1}, \Bar{u}^{k-1}) \nonumber \\
& = \sum_{s\in S_{k}} P^{\gamma}(d\Bar{y}^{k}, s |\omega_{0}, \omega_{s_{0}}, \Bar{y}^{1}, \Bar{u}^{1}, \cdots, \Bar{y}^{k-1}, \Bar{u}^{k-1}) \nonumber \\
& = \sum_{s\in S_{k}} P^{\gamma}\bigg(d\Bar{y}^{k} \bigg| s, \omega_{0}, \omega_{s_{0}}, \Bar{y}^{1}, \Bar{u}^{1}, \cdots, \Bar{y}^{k-1}, \Bar{u}^{k-1}\bigg) P^{\gamma}(s | \omega_{0}, \omega_{s_{0}}, \Bar{y}^{1}, \Bar{u}^{1}, \cdots, \Bar{y}^{k-1}, \Bar{u}^{k-1}) \nonumber \\
& = \sum_{s \in S_{k}} P\bigg(d\Bar{y}^{k} \bigg| s, \omega_{0}, \omega_{s_{0}}, \Bar{y}^{1}, \Bar{u}^{1}, \cdots, \Bar{y}^{k-1}, \Bar{u}^{k-1}\bigg) P(s | \omega_{0}, \omega_{s_{0}}, \Bar{y}^{1}, \Bar{u}^{1}, \cdots, \Bar{y}^{k-1}, \Bar{u}^{k-1}) \label{CUsed} \\
&  =P(d\Bar{y}^{k}|\omega_{0}, \omega_{s_{0}}, \Bar{y}^{1}, \Bar{u}^{1}, \cdots, \Bar{y}^{k-1}, \Bar{u}^{k-1}),
\end{align}
that is, the conditional probability is policy independent. 
Equation (\ref{CUsed}) critically follows from Property C, via (\ref{orderingEqnMap}); as the past control and measurements uniquely determine both the ordering as well as the conditional probability. \qed

\begin{theorem} \label{mainThmNSSR}
Suppose that an N-DM non-sequential stochastic control problem possesses Property C for all admissible policies $\gamma \in \Gamma$, and that the set $\mathbb{Y}^{i}$ is countable for $i = 1, 2, \cdots, N$, then the problem admits a policy-independent static reduction where each DM's measurements are independent from other DMs' policies/actions and are exogenous.
\end{theorem}

\textbf{Proof.}
\newline
{\bf Step 1.}
Let us first write:
\[
P^{\gamma}(d\omega_{0}, d\omega_{s_{0}}, dy, d{\bf u}) = P^{\gamma}(d\omega_{0}, d\omega_{s_{0}})P^{\gamma}(dy, du|\omega_{0}, \omega_{s_{0}})\]
where:
\begin{align}
& P^{\gamma}(dy, du|\omega_{0}, \omega_{s_{0}}) = P^{\gamma}(d\Bar{y}^{1}|\omega_{0}, \omega_{s_{0}})P^{\gamma}(d\Bar{u}^{1}|\Bar{y}^{1}) \nonumber \\
&  \quad \quad \quad \times P^{\gamma}(d\Bar{y}^{2}|\omega_{0}, \omega_{s_{0}}, \Bar{y}^{1}, \Bar{u}^{1})P^{\gamma}(d\Bar{u}^{2}|\omega_{0}, \omega_{s_{0}}, \Bar{y}^{1}, \Bar{u}^{1},\Bar{y}^{2})  \nonumber \\
& \quad \quad \quad \cdots \times P^{\gamma}(d\Bar{y}^{N}|\omega_{0}, \omega_{s_{0}}, \Bar{y}^{1}, \Bar{u}^{1}, \cdots, \Bar{y}^{N-1}, \Bar{u}^{N-1})  \nonumber \\
& \quad \quad \quad  \times P^{\gamma}(d\Bar{u}^{N}| \Bar{y}^{N},\omega_{0}, \omega_{s_{0}}, \Bar{y}^{1}, \Bar{u}^{1}, \cdots, \Bar{y}^{N-1}, \Bar{u}^{N-1})\label{reduc1}
\end{align}

{\bf Step 2.} 
We have that, for all $1 \leq k \leq N$
\[
 P^{\gamma}(d\Bar{y}^{k}|\omega_{0}, \omega_{s_{0}}, \Bar{y}^{1}, \Bar{u}^{1}, \cdots, \Bar{y}^{k-1}, \Bar{u}^{k-1}) = P(d\Bar{y}^{k}|\omega_{0}, \omega_{s_{0}}, \Bar{y}^{1}, \Bar{u}^{1}, \cdots, \Bar{y}^{k-1}, \Bar{u}^{k-1})
 \]
that is, the conditional probability is policy independent. This is due to Lemma \ref{ImaginarySequentialModel}.

{\bf Step 3.} 
We have that, for all $1 \leq k \leq N$
\begin{align}
& P^{\gamma}(d\Bar{u}^{k}| \Bar{y}^{k},\omega_{0}, \omega_{s_{0}}, \Bar{y}^{1}, \Bar{u}^{1}, \cdots, \Bar{y}^{k-1}, \Bar{u}^{k-1})\nonumber \\
&= \sum_{s_k} P^{\gamma}(d\Bar{u}^{k} | s_k , \Bar{y}^{k},\omega_{0}, \omega_{s_{0}}, \Bar{y}^{1}, \Bar{u}^{1}, \cdots, \Bar{y}^{k-1}, \Bar{u}^{k-1}) P^{\gamma}(s_k| \Bar{y}^{k},\omega_{0}, \omega_{s_{0}}, \Bar{y}^{1}, \Bar{u}^{1}, \cdots, \Bar{y}^{k-1}, \Bar{u}^{k-1}) \nonumber \\
&= \sum_{s_k} P^{\gamma}(d\Bar{u}^{k} | s_k , \Bar{y}^{k},\omega_{0}, \omega_{s_{0}}, \Bar{y}^{1}, \Bar{u}^{1}, \cdots, \Bar{y}^{k-1}, \Bar{u}^{k-1}) P(s_k| \Bar{y}^{k},\omega_{0}, \omega_{s_{0}}, \Bar{y}^{1}, \Bar{u}^{1}, \cdots, \Bar{y}^{k-1}, \Bar{u}^{k-1}) \label{useCInd} \\
&= \sum_{s_k} P^{\gamma}(d\Bar{u}^{k} | s_k , \Bar{y}^{k}) P(s_k| \Bar{y}^{k},\omega_{0}, \omega_{s_{0}}, \Bar{y}^{1}, \Bar{u}^{1}, \cdots, \Bar{y}^{k-1}, \Bar{u}^{k-1}) \nonumber \\
&= \sum_{s_k} P^{\gamma^{s_k}}(d\Bar{u}^{k} | \Bar{y}^{k}) P(s_k| \omega_{0}, \omega_{s_{0}}, \Bar{y}^{1}, \Bar{u}^{1}, \cdots, \Bar{y}^{k-1}, \Bar{u}^{k-1}) \label{condIndM2} 
\end{align} 

Here, (\ref{useCInd}) follows also from Property C; note that $s_k$ is determined given the conditioned information. With Property C, the measurements $\Bar{y}^{i}$ are related to the original measurements $y^{i}$ by: 
  \begin{equation}
  \label{II}
  \Bar{y}^{i} = \sum_{k=1}^{N}y^{k}1_{\{\psi_{i}(\omega_{s_{0}}, \Bar{u}^{1}, \cdots, \Bar{u}^{i-1}) = k\}}
  \end{equation}
  and the actions $\Bar{u}^{i}$ are related to $u^{i}$ by:
  \begin{equation}
  \label{IU2}
      \bar{u}^{i} = \sum_{k=1}^{N}u^{k}1_{\{\psi_{i}(\omega_{s_{0}}, \Bar{u}^{1}, \cdots, \Bar{u}^{i-1}) = k\}}
  \end{equation}
Equation (\ref{condIndM2}) follows from Property C.

{\bf Step 4.} Thus, we write (\ref{reduc1}) as
\begin{align}\label{reduc12}
 P^{\gamma}(dy, du|\omega_{0}, \omega_{s_{0}}) &= P(d\Bar{y}^{1}|\omega_{0}, \omega_{s_{0}})P^{\Bar{\gamma}^{1}}(d\Bar{u}^{1}|\omega_{0}, \omega_{s_{0}},\Bar{y}^{1})  \nonumber \\
& \qquad \times P(d\Bar{y}^{2}|\omega_{0}, \omega_{s_{0}},
\Bar{y}^{1}, \Bar{u}^{1})P^{\Bar{\gamma}^{2}}(d\Bar{u}^{2}|\omega_{0}, \omega_{s_{0}}, \Bar{y}^{1}, \Bar{u}^{1},\Bar{y}^{2})  \nonumber \\
& \quad \quad \cdots \times P(d\Bar{y}^{N}|\omega_{0}, \omega_{s_{0}}, \Bar{y}^{1}, \Bar{u}^{1}, \cdots, \Bar{y}^{N-1}, \Bar{u}^{N-1}) \nonumber \\
& \quad \quad \times P^{\Bar{\gamma}^{N}}(d\Bar{u}^{N}| d\Bar{y}^{N},\omega_{0}, \omega_{s_{0}}, \Bar{y}^{1}, \Bar{u}^{1}, \cdots, \Bar{y}^{N-1}, \Bar{u}^{N-1})
\end{align}
We then have that, in view of Step 3,
\begin{align}\label{reduc131313}
& P^{\gamma}(dy, du|\omega_{0}, \omega_{s_{0}}) \nonumber \\
& = P(d\Bar{y}^{1}|\omega_{0}, \omega_{s_{0}}) \bigg( \sum_{s_1} P^{\Bar{\gamma}^{1}}(d\Bar{u}^{1} | \Bar{y}^{1}) P(s_1| \omega_{0}, \omega_{s_{0}})  \bigg)  \nonumber \\
& \qquad \times P(d\Bar{y}^{2}|\omega_{0}, \omega_{s_{0}}, \Bar{y}^{1}, \Bar{u}^{1})\bigg(\sum_{s_2} P^{\Bar{\gamma}^{2}}(d\Bar{u}^{2} | \Bar{y}^{2}) P(s_2| \omega_{0}, \omega_{s_{0}}, \Bar{y}^{1}, \Bar{u}^{1})  \bigg)  \nonumber \\
& \quad \quad \cdots \times P(d\Bar{y}^{N}|\omega_{0}, \omega_{s_{0}}, \Bar{y}^{1}, \Bar{u}^{1}, \cdots, \Bar{y}^{N-1}, \Bar{u}^{N-1}) \nonumber \\
& \quad \quad \times \bigg( \sum_{s_N} P^{\Bar{\gamma}^{N}}(d\Bar{u}^{N} | \Bar{y}^{N}) P(s_N| \omega_{0}, \omega_{s_{0}}, \Bar{y}^{1}, \Bar{u}^{1}, \cdots, \Bar{y}^{N-1}, \Bar{u}^{N-1})   \bigg) 
\end{align}
%
   
The stochastic kernels $P(d\Bar{y}^{i}|\omega_{0}, \omega_{s_{0}}, \Bar{y}^{1}, \cdots, \Bar{y}^{i-1}, \Bar{u}^{1}, \cdots, \Bar{u}^{i-1})$ are policy-independent, as we are conditioning on the actions of all previous DMs as well as the information used to select those actions, thus the policy that maps the past DMs' measurements to their actions is irrelevant. 

 We have that for any $i\in \{1,\cdots,N\}$ one can always write 
\begin{eqnarray*}
    P(ds_{i}|\omega_{0},\omega_{s_{0}},\Bar{y}^{1},\cdots,\Bar{y}^{i-1},\Bar{u}^{1},\cdots,\Bar{u}^{i-1})&=&\int g_{i}(s_{i},\omega_{0},\omega_{s_{0}},\Bar{y}^{1},\cdots,\Bar{y}^{i-1},\Bar{u}^{1},\cdots,\Bar{u}^{i-1})\phi_{i}(ds_{i})
\end{eqnarray*}
Thus, we get that 
\begin{align}
J(\gamma) 
&= \int P(d\omega_{0}, d\omega_{s_{0}})\prod_{i=1}^{N}\phi_{i}(ds_{i})  \nonumber \\
&  \qquad \times f_{1}(\Bar{y}^{1}, \omega_{s_{0}}, \omega_{0})\Bar{Q}_{1}(d\Bar{y}^{1}) P^{\Bar{\gamma}^{1}}(d\Bar{u}^{1} | \Bar{y}^{1}, s_{1}) g_{1}(s_{1}, \omega_{0}, \omega_{s_{0}})    \nonumber \\
& \qquad \times f_{2}(\Bar{y}^{2},\omega_{0}, \omega_{s_{0}}, \Bar{y}^{1}, \Bar{u}^{1}) \Bar{Q}_{2}(d\Bar{y}^{2})
 P^{\Bar{\gamma}^{2}}(d\Bar{u}^{2} | \Bar{y}^{2}, s_{2}) g_{2}(s_{2}, \omega_{0}, \omega_{s_{0}}, \Bar{y}^{1}, \Bar{u}^{1})    \nonumber \\
&\qquad \cdots \qquad \nonumber \\
&  \qquad \times f_{N}(\Bar{y}^{N}, \omega_{s_{0}}, \omega_{0}, \Bar{y}^{1}, \cdots, \Bar{y}^{N-1}, \Bar{u}^{1}, \cdots, \Bar{u}^{N-1})\Bar{Q}_{N}(d\Bar{y}^{N}) \nonumber \\
& \qquad  P^{\Bar{\gamma}^{N}}(d\Bar{u}^{N} | \Bar{y}^{N}, s_{N})g_{N}(s_{N}, \omega_{0}, \omega_{s_{0}}, \Bar{y}^{1}, \Bar{u}^{1}, \cdots, \Bar{y}^{N-1}, \Bar{u}^{N-1})   \nonumber \\
& \qquad \qquad \qquad \qquad \times \tilde{c}(\omega_{0}, \omega_{s_{0}}, s_{[1,N]},\Bar{u}^{[1,N]}),
\end{align}
Now, define $\hat{y}^{i}=(\Bar{y}^{i},s_{i})$ and $\hat{Q}_{i}(d\hat{y}^{i})=\Bar{Q}_{i}(d\Bar{y}^{i})\phi_{i}(ds_{i})$. Then, one can write

\begin{align}
J(\gamma) 
&= \int P(d\omega_{0}, d\omega_{s_{0}})  \nonumber \\
&  \qquad \times f_{1}(\Bar{y}^{1}, \omega_{s_{0}}, \omega_{0})\hat{Q}_{1}(d\hat{y}^{1}) P^{\Bar{\gamma}^{1}}(d\Bar{u}^{1} | \hat{y}^{1}) g_{1}(s_{1}, \omega_{0}, \omega_{s_{0}})    \nonumber \\
& \qquad \times f_{2}(\Bar{y}^{2},\omega_{0}, \omega_{s_{0}}, \Bar{y}^{1}, \Bar{u}^{1}) \hat{Q}_{2}(d\hat{y}^{2})
 P^{\Bar{\gamma}^{2}}(d\Bar{u}^{2} | \hat{y}^{2}) g_{2}(s_{2}, \omega_{0}, \omega_{s_{0}}, \Bar{y}^{1}, \Bar{u}^{1})    \nonumber \\
&\qquad \cdots \qquad \nonumber \\
&  \qquad \times f_{N}(\Bar{y}^{N}, \omega_{s_{0}}, \omega_{0}, \Bar{y}^{1}, \cdots, \Bar{y}^{N-1}, \Bar{u}^{1}, \cdots, \Bar{u}^{N-1})\hat{Q}_{N}(d\hat{y}^{N}) \nonumber \\
& \qquad  P^{\Bar{\gamma}^{N}}(d\Bar{u}^{N} | \hat{y}^{N})g_{N}(s_{N}, \omega_{0}, \omega_{s_{0}}, \Bar{y}^{1}, \Bar{u}^{1}, \cdots, \Bar{y}^{N-1}, \Bar{u}^{N-1})   \nonumber \\
& \qquad \qquad \qquad \qquad \times \tilde{c}(\omega_{0}, \omega_{s_{0}}, s_{[1,N]},\Bar{u}^{[1,N]}),
\end{align}
The last expression can be simplified as 
\begin{align}
J(\gamma) 
&= \int P(d\omega_{0}, d\omega_{s_{0}})   \times \hat{Q}_{1}(d\hat{y}^{1}) P^{\Bar{\gamma}^{1}}(d\Bar{u}^{1} | \hat{y}^{1})    \times  \hat{Q}_{2}(d\hat{y}^{2})
 P^{\Bar{\gamma}^{2}}(d\Bar{u}^{2} | \hat{y}^{2})   \cdots  \nonumber \\
&  \qquad \times \hat{Q}_{N}(d\hat{y}^{N})   P^{\Bar{\gamma}^{N}}(d\Bar{u}^{N} | \hat{y}^{N})    \times {\mathbf{c}}(\omega_{0}, \omega_{s_{0}}, \hat{y}^{[1,N]},\Bar{u}^{[1,N]}), \label{rep}
\end{align}

where,

\begin{eqnarray*}
   && \mathbf{c}(\omega_{0}, \omega_{s_{0}}, \hat{y}^{[1,N]},\Bar{u}^{[1,N]})= f_{1}(\Bar{y}^{1}, \omega_{s_{0}}, \omega_{0})\times g_{1}(s_{1}, \omega_{0}, \omega_{s_{0}}) \times \cdots \\
   &&\qquad \times f_{N}(\Bar{y}^{N}, \omega_{s_{0}}, \omega_{0}, \Bar{y}^{1}, \cdots, \Bar{y}^{N-1}, \Bar{u}^{1}, \cdots, \Bar{u}^{N-1})\times  g_{N}(s_{N}, \omega_{0}, \omega_{s_{0}}, \Bar{y}^{1}, \Bar{u}^{1}, \cdots, \Bar{y}^{N-1}, \Bar{u}^{N-1}) \\
   && \qquad \qquad \qquad \times  \tilde{c}(\omega_{0}, \omega_{s_{0}}, s_{[1,N]},\Bar{u}^{[1,N]})
\end{eqnarray*}

Thus, we have an equivalent static problem with independent measurements (generated under measures $\hat{Q}_{1}(d\Bar{y}^{1}), \cdots,  \hat{Q}_{N}(d\Bar{y}^{N})$), and the formulation is policy-independent. 
\qed

We note that in the representation (\ref{rep}), each DM's measurement is exogenous and independent of other DMs' policies. This brings the problem back to its original measurement/policy domain with original measurement spaces.
\begin{remark}
Note that when the cost function $c$ is permutation invariant, with
\[c(\omega_{0}, \omega_{s_{0}},  u_1,\cdots,u_N) =: c(\omega_{0}, \sigma^{-1}_{s}(\Bar{u}^{[1,N]}))\]
the reduction simplifies, as we will see in the example in Section \ref{exampleSecS}.
\end{remark}

\begin{remark}
One further method for converting non-sequential problems into sequential ones was proposed by Andersland in \cite{andersland1991decoupling}, where again policies appear as cost parameters; provided that the problem satisfies Property CI.  In this approach, a new agent, acting before all others, is introduced which simulates the actions of the future agents whose incompatible actions are given a zero reward. The simulation effectively declares the team policy. 
\end{remark}

\section{Example}\label{exampleSecS}

We present in this section an example to illustrate the non-sequential policy-independent static reduction process. In particular, we revisit the example given in Example \ref{Ornek1}. 

With the information structure as given in Example \ref{Ornek1}, define the cost function as: $c(\omega_{0}, {\bf u}) = \omega_{0} + u^{1} + u^{2} + u^{3}$. Then, let us write:
    
    $$P(\Bar{y}^{1} \in A|\omega_{s_{0}}, \omega_{0}) = \int_{A}f_{1}(\Bar{y}^{1}, \omega_{s_{0}}, \omega_{0})\Bar{Q}_{i}(d\Bar{y}^{i})= \sum_{\Bar{y}^{1} \in A} \frac{P(\Bar{y}^{1}, \omega_{s_{0}}, \omega_{0})}{P(\omega_{s_{0}}, \omega_{0})}$$
    $$= \int_{A} \frac{P(\Bar{y}^{1}|\omega_{s_{0}}, \omega_{0})P(\omega_{s_{0}}, \omega_{0})}{P(\omega_{s_{0}}, \omega_{0})\mu^{1}(\Bar{y}^{1})}\mu^{1}(d\Bar{y}^{1})$$
    
    \noindent
    Where $\mu^{1}$ is the uniform distribution on $\Bar{Y}^{i}$, and $P(\Bar{y}^{1} \in A|\omega_{s_{0}}, \omega_{0})$ is a well defined stochastic kernel representing either the stochastic mapping $g_{1}$ or $g_{2}$. For example, if $\omega_{s_{0}} = 1$ and $\omega_{0} = 1$, then $P(\Bar{y}^{1} = 1|\omega_{s_{0}}, \omega_{0}) = P(y^{2} = 1|\omega_{s_{0}}, \omega_{0}) = 1$. That is, the mapping is deterministic since $\omega_{2}$ is irrelevant in this case. More interestingly, in the case that $\omega_{s_{0}} = 1$ and $u^{2} = 1$, we have that:
    
    $$P(\Bar{y}^{2} = 0| \omega_{s_{0}}, \omega_{0}, \Bar{y}^{1}, \Bar{u}^{1}) = P(y^{1} = 0| \omega_{s_{0}}, \omega_{0}, y^{2}, u^{2}) = P(\omega_{1} = 1)$$
Now, we can construct the change in measure formula similarly for $\Bar{y}^{2}$ and $\Bar{y}^{3}$: 
    $$P(\Bar{y}^{2} \in A|\omega_{s_{0}}, \omega_{0}, \Bar{y}^{1}, \Bar{u}^{1})=\int_{A} \frac{P(\Bar{y}^{2}, \omega_{s_{0}}, \omega_{0}, \Bar{y}^{1}, \Bar{u}^{2})}{P(\omega_{s_{0}}, \omega_{0}, \Bar{y}^{1}, \Bar{u}^{1})\mu^{2}(\Bar{y}^{2})}\mu^{2}(d\Bar{y}^{2})$$
And:

$$P(\Bar{y}^{3} \in A|\omega_{s_{0}}, \omega_{0}, \Bar{y}^{1}, \Bar{y}^{2}, \Bar{u}^{1}, \Bar{u}^{2}) = \int_{A} \frac{P(\Bar{y}^{3}, \omega_{s_{0}}, \omega_{0}, \Bar{y}^{1}, \Bar{y}^{2}, \Bar{u}^{1}, \Bar{u}^{2})}{P(\omega_{s_{0}}, \omega_{0}, \Bar{y}^{1}, \Bar{y}^{2}, \Bar{u}^{1}, \Bar{u}^{2})\mu^{2}(\Bar{y}^{3})}\mu^{2}(d\Bar{y}^{3})$$

\noindent
So now the total cost can be written as:

$$J(\gamma) = \int P(d\omega_{0}, d\omega_{s_{0}})\frac{P(\Bar{y}^{1}, \omega_{s_{0}}, \omega_{0})}{P(\omega_{s_{0}}, \omega_{0})\mu^{1}(\Bar{y}^{1})}\mu^{1}(d\Bar{y}^{1})\frac{P(\Bar{y}^{2}, \omega_{s_{0}}, \omega_{0}, \Bar{y}^{1}, \Bar{u}^{2})}{P(\omega_{s_{0}}, \omega_{0}, \Bar{y}^{1}, \Bar{u}^{1})\mu^{2}(\Bar{y}^{2})}\mu^{2}(d\Bar{y}^{2})$$
$$\times \frac{P(\Bar{y}^{3}, \omega_{s_{0}}, \omega_{0}, \Bar{y}^{1}, \Bar{y}^{2}, \Bar{u}^{1}, \Bar{u}^{2})}{P(\omega_{s_{0}}, \omega_{0}, \Bar{y}^{1}, \Bar{y}^{2}, \Bar{u}^{1}, \Bar{u}^{2})\mu^{3}(\Bar{y}^{3})}\mu^{3}(d\Bar{y}^{3})(\omega_{0}+u^{1}+u^{2}+u^{3})$$
Where $u^{i} = \Bar{u}^{1}1_{\{\psi_{1}(\omega_{s_{0}}) = i\}} + \Bar{u}^{2}1_{\{\psi_{2}(\omega_{s_{0}}, \Bar{u}^{1}) = i\}} + \Bar{u}^{3}1_{\{\psi(\omega_{s_{0}}, \Bar{u}^{1}, \Bar{u}^{2}) = i\}}$, and since each $\Bar{u}^{i}$ appears in the cost function with equal contribution (and the cost is permutation invariant), we can write:

$$J^{*}(\gamma^{*}) = \int P(d\omega_{0}, d\omega_{s_{0}})\frac{P(\Bar{y}^{1}, \omega_{s_{0}}, \omega_{0})}{P(\omega_{s_{0}}, \omega_{0})\mu^{1}(\Bar{y}^{1})}\mu^{1}(d\Bar{y}^{1})\frac{P(\Bar{y}^{2}, \omega_{s_{0}}, \omega_{0}, \Bar{y}^{1}, \Bar{u}^{2})}{P(\omega_{s_{0}}, \omega_{0}, \Bar{y}^{1}, \Bar{u}^{1})\mu^{2}(\Bar{y}^{2})}\mu^{2}(d\Bar{y}^{2})$$
$$\times \frac{P(\Bar{y}^{3}, \omega_{s_{0}}, \omega_{0}, \Bar{y}^{1}, \Bar{y}^{2}, \Bar{u}^{1}, \Bar{u}^{2})}{P(\omega_{s_{0}}, \omega_{0}, \Bar{y}^{1}, \Bar{y}^{2}, \Bar{u}^{1}, \Bar{u}^{2})\mu^{3}(\Bar{y}^{3})}\mu^{3}(d\Bar{y}^{3})(\omega_{0}+\Bar{u}^{1}+\Bar{u}^{2}+\Bar{u}^{3})$$
Now, the new cost function, $c_{s}$ is defined as:

\begin{align} \nonumber c_{s}(\omega_{0}, \omega_{s_{0}}, y, {\bf u}) = \frac{P(\Bar{y}^{1}, \omega_{s_{0}}, \omega_{0})}{P(\omega_{s_{0}}, \omega_{0})\mu^{1}(\Bar{y}^{1})}\frac{P(\Bar{y}^{2}, \omega_{s_{0}}, \omega_{0}, \Bar{y}^{1}, \Bar{u}^{2})}{P(\omega_{s_{0}}, \omega_{0}, \Bar{y}^{1}, \Bar{u}^{1})\mu^{2}(\Bar{y}^{2})} \\ \nonumber \times \frac{P(\Bar{y}^{3}, \omega_{s_{0}}, \omega_{0}, \Bar{y}^{1}, \Bar{y}^{2}, \Bar{u}^{1}, \Bar{u}^{2})}{P(\omega_{s_{0}}, \omega_{0}, \Bar{y}^{1}, \Bar{y}^{2}, \Bar{u}^{1}, \Bar{u}^{2})\mu^{3}(\Bar{y}^{3})}(\omega_{0}+\Bar{u}^{1}+\Bar{u}^{2}+\Bar{u}^{3}),\end{align}

\noindent
and the new joint measure on $(d\omega_{0}, d\omega_{s_{0}}, dy)$ is $P(d\omega_{0}, d\omega_{s_{0}})\mu^{1}(d
\Bar{y}^{1})\mu^{2}(d\Bar{y}^{2})\mu^{3}(\Bar{y}^{3})$. We could further simplify the reduction as in Step 5 in the proof of Theorem \ref{mainThmNSSR}.

\section{Conclusion} In this paper, we revisited several seminal studies on non-sequential decentralized stochastic control. We studied the notion of Causality, presented an equivalent formulation, and showed that Causality is equivalent to Causal Implementability (and Dead-Lock Freeness) for standard Borel models. We showed that Causality, under an absolute continuity condition, allows for an equivalent static model whose reduction is policy-independent. We further showed that under more relaxed conditions on the model, such as solvability, such a reduction (when possible) is policy-dependent or includes policies as parameters in the cost of the reduced model, and thus has limited utility. 




\bibliographystyle{siam}

\end{document}